\documentclass[preprint, 10pt]{elsarticle}
\usepackage[centertags]{amsmath}
\usepackage{amssymb}
\usepackage{amsthm}
\usepackage{amsfonts,amssymb}
\usepackage{enumerate}
\usepackage{mathrsfs}
\usepackage{graphicx}

\journal{}

\theoremstyle{plain} 

\newtheorem{thm}{Theorem}

\theoremstyle{definition}

\theoremstyle{remark}
  
\newtheorem*{case1.1}{\textbf{\small Case 1.1}}

\DeclareMathOperator{\sech}{sech}

\DeclareMathOperator{\pd}{\partial}

\DeclareFontFamily{OT1}{pzc}{}
\DeclareFontShape{OT1}{pzc}{m}{it}{<-> s * [1.10] pzcmi7t}{}
\DeclareMathAlphabet{\mathpzc}{OT1}{pzc}{m}{it}

\DeclareMathSymbol{\R}{\mathalpha}{AMSb}{"52}
\DeclareMathSymbol{\C}{\mathalpha}{AMSb}{"43}

\newcommand{\eps}{\varepsilon}
\newcommand{\set}[1]{\left\{#1\right\}}

\newcommand{\comment}[1]{}

\newcommand{\bv}{\mathbf{v}}

\newcommand{\bsube}{\begin{subequations}}
\newcommand{\esube}{\end{subequations}}

\begin{document}

\begin{frontmatter}



\title{Group classification and exact solutions of a class of nonlinear waves}

\author{J.C. Ndogmo\corref{cor1}}
\ead{jean-claude.ndogmo@univen.ac.za}

\address{Department of Mathematical and computational sciences\\
University of Venda\\
P/B X5050, Thohoyandou 0950, South Africa}

\cortext[cor1]{Corresponding author}


\begin{abstract}
We apply an extension of a new method of group classification to a family of nonlinear wave equations labelled by two arbitrary functions, each depending on its own argument. The results obtained confirm the efficiency of the proposed method  for group classification, termed the method of indeterminates. A model equation from the classified  family of fourth order Lagrange equations is singled out. Travelling wave solutions of the latter are found through a similarity reduction by variational symmetry operators, followed by a double order reduction into a second order ordinary differential equation. Multi-soliton solutions and other exact solutions  are also found by various methods including Lie group and  Hirota methods. The most general action of the full symmetry group on any given solution is provided.  Some remarkable  facts on Lagrange equations emerging from the whole study are outlined.
\end{abstract}

\begin{keyword}
Nonlinear waves \sep Group classification \sep Travelling waves \sep Multi-solitons \sep Symmetry properties

\MSC[2020] 58D19 \sep  76M60 \sep 35C05   \sep 35C08
\end{keyword}

\end{frontmatter}

\section{Introduction}
\label{s:intro}

The modelling of physical phenomena by differential equations is usually based on parameter approximations or on very specific and simplified cases of the phenomenon. Such a modelling therefore yields results expressed in terms of constant parameters or more generally numeric values, which in reality could be replaced by more global expressions, including often arbitrary functions that may depend on all the variables, both the dependent and independent ones involved in the equation. The group classification of such families of equations depending on arbitrary parameters can  be used to completely determine the symmetry properties of each symmetry class from the equation. This yields in particular many of the relevant information on the whole family of equations, some of which might not be apparent in the simplified models. In the context of Lie group analysis, such information will usually include amongst others symmetry algebras,  various types of group-invariant solutions  and their properties, equivalence groups, conservation laws and conserved quantities, as well as symmetry integrability properties.\par

   Implementing the group classification of differential equations has proved however to be quite challenging, first because the direct method of analysis usually employed \cite{ovsyC,gungorC,ndogC,ZdanovC} involves lengthy case analyses and is thus prone to incorrect  results \cite{bagderC,popovRDE,ndog_hyperb}. Also, there is no known methods up to now to predict the approximate number of symmetry classes contained in a given family of parameterized equations, and even under equivalence transformations the implementation of a group classification can turn out to be far more complicated and lengthy than anticipated. \par

Although the group classification of differential equations has usually been done through an improvised method so-called the direct method, other methods have been proposed and used, such as the algebraic method \cite{gungorC,moyoC2} which is quite popular.  Another classification method that has emerged in recent years, often called furcate splitting,  has been applied in some papers (see  \cite{popovRDE,oovaneeva} and the references therein). Another group classification method that has appeared more recently is probably the one coined the method of indeterminates, and it has been applied to the group classification of a family of hyperbolic equations \cite{ndog_hyperb,classifIOP21}.  However, in the simpler case of the study carried out in \cite{ndog_hyperb}, the family of equations was labelled by a single arbitrary function of a single variable. On the other hand,  the application of the method carried out in \cite{classifIOP21} was very brief. \par

The present paper aims to extend the application of our method of indeterminates for  group classification to a more general family of equations, namely  one containing two arbitrary functions, each of which depends on a different variable. More precisely, we consider in this paper the problem of group classification of a family of nonlinear wave equations of the form

\begin{equation}\label{e:main}
u_{tt}  - (A\, u_{xxxx} + S' u_x u_{xx}  + 2 S u_{xx} + Z)=0,
\end{equation}
where $u=u(t,x).$ Moreover, $u_x= \pd u/ \pd x,\, u_y= \pd u/ \pd y, \, u_{xy}= \pd^2 u/ (\pd x\, \pd y)$ and so on. Similarly, a prime on a function of one variable denotes the function's derivative  with respect to its argument. Additionally,  $Z=Z(u)$ and  $S=S(u_x)$ are arbitrary functions while $A$ is an arbitrary nonzero constant. The summary of this classification result is presented in Table \ref{tb:1}.\par

In order to gain more insight into the symmetry as well as the physical properties inherent to the family of equations \eqref{e:main}, a model equation from the list of symmetry classes obtained is singled out for further studies, namely equation \eqref{e:5th} in Section \ref{s:soln}.
    Given that the latter equation is Lagrangian, we obtain a double order reduction of its similarity reduced version as well as  expressions for the related group-invariant solutions. Moreover, given  that \eqref{e:5th}, and in fact the whole of \eqref{e:main},  is invariant under both time and space translations, explicit travelling wave solutions which are actually 1-solitons, as well as other explicit solutions are also constructed. Independently from all the above methods, families of one- and two-soliton solutions are derived by the Hirota method. \par

    Similarly, to gain more insight into the properties of the solutions of the model equation thus obtained an expression is given for the most general action of the full symmetry group on a given solution. This Lie group action solves amongst others the crucial problem of genuinely generating new soliton solutions from given ones. Indeed, in the current literature no true generation of solitions from given ones seems to have been performed and the expression 'generation of soliton solutions' usually refers only to the complete determination from scratch  of  solutions (see e.g. \cite{demontis,yao} and the references therein).\par

   In the quest for more discoveries, the family of equations \eqref{e:main} was chosen to be not only of the fourth order, but also of Lagrange type, as both types of equations are seldom studied, especially those of Euler-Lagrange type. This choice has however led to some remarkable facts emerging from the whole study concerning Lagrange equations and outlined  in the concluding section of the paper.


\section{Equivalence group } \label{s:equivGrp}

The equivalence group of \eqref{e:main} will be the Lie pseudo-group of point transformations that leaves \eqref{e:main} invariant up to arbitrary parameter functions. That is, it leaves the family of equations \eqref{e:main} invariant. To a given monomial in dependent variables such as $u$ together with its derivatives, we assign a number call its weight, and which equals the total order of derivatives of $u$ in the given monomial. In this way,  $u,\, u_x^2 u_{xx}u,\, u_{xt}u_t^3$ have weight zero, four and five, respectively.  Given that a point transformation maps a monomial of maximal weight in the original equation to a monomial of maximal weight in the transformed equation, equivalence transformations should preserve the maximal weight. Based on these considerations, it is easily seen that  equivalence transformations of \eqref{e:main} must satisfy
\begin{equation}\label{e:Gequiv}
u= H(y,z, w),\; x= G(y,z),\;  t= F(y,z), \;  \text{ where } (F_z G_y- F_y G_z)H_w \neq 0,
\end{equation}
and  $w=w(y,z),$ and for some functions $F=F(y,z), G=G(y,z),$ and $H= H(y,z)$ to be specified. For \eqref{e:Gequiv} to be indeed an equivalence transformation it is in particular necessary that the coefficients of every new monomial in $w$ and its derivatives occurring in the transformed version of \eqref{e:main} under \eqref{e:Gequiv} vanish. For instance, the vanishing of the coefficient $A F_z^4(F_z G_y- F_y G_z)^3 H_w$ of $w_{yyyy}$ in the transformed version of \eqref{e:main} shows that $F=F(y)$ should depend only on $y.$ Proceeding in this way for all other new monomials in the transformed equation, and more generally, requesting that the transformed version of \eqref{e:main} be of the same form as \eqref{e:main} yields the following result.
\begin{thm}
The equivalence group of \eqref{e:main} consists of invertible linear transformations
\begin{subequations}\label{e:equivM}
\begin{align}
t &= J z + a,\quad x= K y + b,\quad u= T w + \gamma, \label{e:equivM1} \\
\intertext{where $J, K, T, \text{ and } a, b, \gamma$ are constants with $0\neq J\,K\,T.$ The transformed version of \eqref{e:main} under \eqref{e:equivM1} takes the form}
w_{yy}&= \frac{A J^2}{K^4} w_{zzzz} + S_1'(w_z) w_z w_{zz} + 2 S_1(w_z) w_{zz}+  Z_1(w) \label{e:equivM2} \\
\intertext{where}
S_1(w_z) &=\frac{J^2}{K^2} S(T w_z/ K),\quad  \text{ and } Z_1(w) =  \frac{J^2}{T} Z(T w+ \gamma).  \label{e:equivM3}
\end{align}
\end{subequations}
\end{thm}

It is clear that the transformations \eqref{e:equivM1} induce on the set of all pairs of functions $S(u_x)$ and $Z(u)$ an equivalence relation which we shall denote  by $\sim\, ,$  and under which any given such pair is equivalent to the pair $S_1(w_z)$ and $ Z_1(w)$ given by \eqref{e:equivM3}, or more simply after a renaming of variables, to the pair $(S_1(u_x), Z_1(u)).$  In particular we can now write $(S, Z) \sim (S_1, Z_1),$ and more loosely $S\sim S_1$ and  $Z\sim Z_1.$ \par

Since the transformations \eqref{e:equivM1} act simultaneously on any given pair $(S_1, Z_1)$  of functions parameters of \eqref{e:main}, in principle the transformation of any one of the components in $(S_1, Z_1),$ should limit the degree of freedom in the transformation of the other function. However, due to the fact that there is a sufficiently large number of free parameters in the expressions \eqref{e:equivM3}, each of the functions $Z_1$ and $S_1$ may be fully transformed independently of the other, as if  \eqref{e:equivM1} were actually not acting simultaneously on pairs of functions.


\section{Lie group classification} \label{s:classif}
Let us recall that the group classification consists in finding all non equivalent symmetry classes according to the values assumed by the labelling parameters of the family of equations. According to the classical method of Lie, this classification can be achieved by an analysis of the determining equations \cite{olv93}. Let us denote by
\begin{equation}\label{e:Gsym}
\bv= \xi \pd_t + \eta \pd_x + \phi \pd_u
\end{equation}
an arbitrary vector in the symmetry algebra of \eqref{e:main}, in which $\xi, \eta$ and $\phi$ each represent functions of $t,x$ and $u.$  Let $\bv^{[4]}$ denote the prolongation of order four of $\bv.$ Then the expression

\begin{equation}\label{e:invSurf1}
\bv^{[4]} (\Delta)_{\big\vert_{\,\Delta=0}} =0,
\end{equation}
where $\Delta$ represents the left hand side of \eqref{e:main}, expanded into a polynomial in the unconstrained derivatives of $u=u(t,x)$ gives rise to so-called determining equations for the components $\xi, \eta,$ and $\phi$ of $\bv,$  and hence determines $\bv.$ In the actual case of \eqref{e:main}, the only constrained derivative of $u$ is $u_x.$ Solving in \eqref{e:invSurf1} all equations not involving  the parameter functions $Z$ and $S$ and their derivatives shows that the functions  $\xi, \eta,$ and $\phi$ should satisfy
\begin{align}
\xi &=  2 k_1 t + k_2,\qquad \eta= k_1 x + k_3,\qquad \phi= k_4 u + h, \label{e:deteqM0}
\end{align}
where here and in the sequel, the $k_j \text{ (for $j=1,\dots,4$)}$ are arbitrary constants and $h=h(t,x).$ The remaining equations of \eqref{e:invSurf1} are then reduced to
\begin{subequations}\label{e:deteqM1}
\begin{align}
0 &= -4 k_1 S -3 h_x S' + (k_1 -k_4)u_x^2 S'' + u_x \left( (k_1 - 3 k_4) S' - h_x S''\right)  \label{e:deteqM1a} \\
0 &= (- 4 k_1 + k_4) Z - h Z' - k_4 u Z' + h_{tt} - 2 S h_{xx} - S' h_{xx} - A h_{xxxx}\label{e:deteqM1b}
\end{align}
\end{subequations}
and are referred to as the classifying equations. Assuming that $Z$ and $S$ are arbitrary in \eqref{e:deteqM1} gives $0=k_1=k_4=h.$ In other words, the principal algebra $L_p$ of \eqref{e:main} is spanned by $\pd_t$ and $\pd_x,$ revealing the fact that as a wave equation \eqref{e:main} admits travelling wave solutions. The group classification of \eqref{e:main} will consist in finding all possible extensions of $L_p$ according to the values of the parameter functions $Z$ and $S.$ However, in  terms of the components $\xi, \eta,$ and $\phi$ as given by   \eqref{e:deteqM0}, this group classification will simply consist in finding all possible values of $h=h(t,x)$ and confirming those of the constants $k_1, \dots, k_4$ for which the  corresponding vector field $\bv$ does yield a symmetry generator.


The group classification of \eqref{e:main} will be performed using an extended version of the method we recently proposed in \cite{ndog_hyperb}, where it was applied to a class of hyperbolic equations in which the classifying equations consisted, amongst others of only one scalar equation.

Roughly speaking, the underlying idea of the method of indeterminates is to consider the classifying equations wherever feasible as a system of polynomial equations in any variables that may be appropriately considered as indeterminates. The method will consist in the actual case in viewing the classifying equations as a polynomial with the monomials having as indeterminates the variables $u$ and $u_x,$ as well as any of the functions $Z=Z(u),$ $S=S(u_x)$ and their derivatives occurring in the classifying equations. We then exploit the fact that any extension of the principal algebra occurs only if one or more of the resulting monomials are linearly dependent.\par

The assumption that a set of $p$ monomials in a given equation is linearly dependent will give rise to a specification of the arbitrary functions, up to some arbitrary constants. Such a specification will be achieved through the solving of a low-order linear differential equation, namely one of order at most two in the actual case, and the functions thus specified will be termed admissible. The resulting  arbitrary constants just mentioned might also be all or partly eliminated through equivalence transformations. The integer $p$ will run from one to $m,$ where $m$ is the maximal number of monomials in the classifying equations.\par

Each resulting admissible function will be substituted into the entire system of equations,  to determine particular values of  arbitrary constants, if any, and to recursively apply the same process for finding admissible functions among any of the yet unspecified arbitrary functions and their subsequent substitution into the entire system of equations. For each possible combination of admissible functions, each of which corresponds to particular values  of each of the arbitrary functions, the corresponding symmetry algebra will then be computed to find out if it does extend the principal algebra.\par

   In the case of the classifying equations \eqref{e:deteqM1}, it is appropriate to start the search of admissible functions with the first equation \eqref{e:deteqM1a}, as it involves only a single parameter function and all possible indeterminates involving this parameter from the entire system of equations. Letting $p=2,$ the vanishing of all possible linear combinations of $p$ monomials yields  after an application of the equivalence transformations \eqref{e:equivM1} the following adjusted list of admissible functions for $S= S(u_x).$

\begin{subequations}\label{e:admFct0S}
\begin{align}
S &= \beta,\; u_x + \beta,\; u_x^\alpha + \beta,\;  u_x^\sigma + \rho\,  u_x^\tau,\; \text{ or } \rho\,  e^{u_x}, \label{e:admFct1S} \\
\intertext{where the arbitrary constants $\alpha, \beta, \rho, \sigma, \text{ and } \tau$ satify}
\alpha &\neq 0, 1;\quad \sigma + \tau =1;\quad  \sigma\, \tau \neq 0;\quad \rho \neq 0. \label{e:admFct2S}
\end{align}
\end{subequations}

The substitution turn by turn in \eqref{e:deteqM1a}  of each particular value of $S$ occurring  in \eqref{e:admFct1S} will determine any particular value of the associated arbitrary constants, and by the same process as for $S$ all admissible functions of $Z= Z(u)$ resulting from the assumption of linear dependence of the monomials in the resulting equation \eqref{e:deteqM1b}. This process will give rise to combinations of pairs $(S, Z)$ of admissible functions, each of which is used to test the extension of the principal algebra.

For an explicit example of application of this procedure, let us consider the case where the admissible value assumed by $S$  from the list \eqref{e:admFct1S} is given by $S= u_x + \beta.$ This reduces \eqref{e:deteqM1} to
\begin{subequations}\label{e:deteqEx1}
\begin{align}
  0 &= -4 k_1 \beta -3 h_x - 3(k_1+ k_4)u_x \label{e:deteqEx1a}\\
   0 &=  (-4 k_1 + k_4) Z - h Z' - k4 u Z' + h_{tt} - 2 \beta h_{xx} - 3 u_x h_{xx} -A\,  h_{xxxx}. \label{e:deteqEx1b}
\end{align}
\end{subequations}
It then follows from \eqref{e:deteqEx1a} that
\[
k_1= -k_4,\qquad \text{ and } h(t,x) = \frac{4}{3}  k_4 x \beta + f,
\]
for some function $f= f(t).$ Substituting the latter values for $k_1$ an $h$ into \eqref{e:deteqEx1} yields
\begin{align}
   0 &= 5 k_4 Z - k_4 u Z' + (-f -  \frac{4}{3}  k_4 x \beta ) Z' + f_{tt}. \label{e:deteqEx2}
\end{align}
To pursue this analysis, we first need to determine the admissible values of the function $Z=Z(u)$ in \eqref{e:deteqEx2}. Given the monomials $u, Z, u Z', Z'$  and $1$ of the polynomial appearing in the right hand side of \eqref{e:deteqEx2} and the equivalence transformations \eqref{e:equivM1}, the admissible values of $Z$ for all possible $p$-tuples of linearly dependent monomials are given in this case by
\begin{equation}\label{e:admZex1}
Z= u,\quad \delta,\quad \log(u),\quad e^u,\quad \text{ or }\quad  u^{\theta},
\end{equation}
where $\delta$ and $\theta$ are arbitrary constants with $\theta \neq 0,1.$
Each value of $Z$ from the list in \eqref{e:admZex1} is to be substituted turn by turn into \eqref{e:deteqEx2} to check for any possible extension of the principal algebra. Letting for instance $Z= e^\theta$ reduces  \eqref{e:deteqEx2} to
\begin{equation}\label{e:deteqEx3}
 k_4 u^\theta (\theta-5) - \frac{1}{3} u^{\theta-1} (3 f + 4 k_4 x \beta) \theta + f'' =0.
\end{equation}
The latter equation shows that any extension of $L_p$ occurs only if $\theta=5$ and $\beta=0,$ and the corresponding symmetry algebra for the resulting pair $(S,Z)= (u_x, u^5)$ has generators
\[
\bv_1 = \pd_t,\quad \bv_2= \pd_x,\quad \bv_3 = 2t \pd_t+ x \pd_x -u \pd_u.
\]
Proceeding in this way for all combinations of pairs of admissible functions $(S, Z)$ yields the group classification of \eqref{e:main}, listed in Table \ref{tb:1}.\par

\comment{
Let us stress that our classification of \eqref{e:main} is meant to be only illustrative of our proposed method rather to be an exhaustive one.  Indeed,  admissible functions $S$ listed in \eqref{e:admFct0S} do not include those from the long list of special functions depending on several free parameters that also occur as admissible functions, and which also include as particular cases a large number of non-special functions. An ample discussion of admissible special functions $S$ appears in the appendix. It is noteworthy that none of the admissible  values of $Z$ occurs as a special function, contrary to case of the parameter function $S= S(u_x).$
}

\begin{table}[hbt!]
\caption{ \label{tb:1}  \protect {\footnotesize \bf  Symmetry classes for Equation \eqref{e:main}.} \footnotesize The  $L_j,\; (j=0, 1, 2),$ each denotes the set of symmetry generators for the  classes  $1$ to $3,$ respectively. In particular, $L_0$ generates the principal algebra. The set of generators for any other  symmetry class is given as an extension of either $L_0, L_1,$ or $L_2$ and represented in the table by simply juxtaposing sequentially the additional symmetry generator to the corresponding generating set $L_j,$ or by  the general expression of a symmetry generator as in the case of Class 8.}

\begin{minipage}[l]{\textwidth}
%

\begin{tabular}{r c c l} \hline \\[-1.5mm]
{\scriptsize \bf No }& {\scriptsize  \bf $S$} & {\scriptsize  \bf $Z$} & {\scriptsize  \textbf{Symmetry generators}}\\[1.5 pt]  \hline

{\scriptsize   1.}& {\scriptsize $S(u_x)$ } & {\scriptsize $Z(u)$  } & {\scriptsize  $\set{\pd_t,\, \pd_x}= L_0 $ } \\

{\scriptsize   2.}& {\scriptsize $S(u_x)$ } & {\scriptsize $\delta$  } & {\scriptsize  $\set{\pd_t,\, \pd_x,\, \pd_u,\, t \pd_u} = L_1  $ } \\

{\scriptsize   3.}& {\scriptsize $S(u_x)$ } & {\scriptsize $u$  } & {\scriptsize $\set{\pd_t,\, \pd_x,\, e^t \pd_u,\,  e^{-t}\pd_u }= L_2$} \\

{\scriptsize   4.}& {\scriptsize $u_x + \beta$ } & {\scriptsize $\theta$  } & {\scriptsize  $L_1, \frac{4 t}{5 \theta } \pd_t + \frac{2x}{5 \theta }\pd_x + \frac{15 t^2 \theta -2 (4 x \beta + 3 u)}{15 \theta } \pd_u $ } \\

{\scriptsize  5.}& {\scriptsize  $u_x + \beta$} & {\scriptsize  $0$ } & {\scriptsize $L_1, -\frac{3 t}{2 \beta } \pd_t - \frac{3x}{4 \beta  }\pd_x + \left(x+  \frac{3 u}{4 \beta  }\right) \pd_u $ }    \\

{\scriptsize  6.}& {\scriptsize   $u_x$} & {\scriptsize  $0$ } & {\scriptsize $L_1,\; 2t \pd_t + x \pd_x -u \pd_u$   } \\

{\scriptsize  7.}& {\scriptsize  $u_x$ } & {\scriptsize $u^5$ } & {\scriptsize $L_0,\; 2t \pd_t + x \pd_x -u \pd_u$     } \\

{\scriptsize  8.}& {\scriptsize  $\beta + \beta_0 u_x^{-2}$ } & {\scriptsize  $\delta + \delta_0 u$ } &
{\scriptsize $k_1 \pd_t  + k_2 \pd_x + (F + k_3 u)\pd_u,\;$   } \\

{\scriptsize  9.}& {\scriptsize $ \beta_0 u_x^{-2}$  } & {\scriptsize  $e^u$ } & {\scriptsize $L_0,\;   2t \pd_t + x \pd_x -4 \pd_u$   } \\

{\scriptsize  10.}& {\scriptsize $ \beta_0 u_x^{-2}$  } & {\scriptsize $u^{\alpha}$  } & {\scriptsize $L_0,\;   2t \pd_t + x \pd_x -\frac{4u}{\alpha-1} \pd_u$   } \\

{\scriptsize  11.}& {\scriptsize   $u_x^\alpha$ } & {\scriptsize  $0$ } & {\scriptsize $L_1,\;   2t \pd_t + x \pd_x +\frac{(\alpha-2)u}{\alpha} \pd_u$    } \\

{\scriptsize  12.}& {\scriptsize  $u_x^{-2/3}$ } & {\scriptsize $\delta$ } & {\scriptsize  $L_1,\;   2t \pd_t + x \pd_x +4 u \pd_u$    } \\

{\scriptsize  13.}& {\scriptsize  $u_x^2$ } & {\scriptsize $\delta$ } & {\scriptsize $L_1,\;   \frac{t}{\delta} \pd_t + \frac{x}{2\delta} \pd_x + t^2 \pd_u$ } \\

{\scriptsize  14.}& {\scriptsize  $u_x^2$  } & {\scriptsize $0$  } & {\scriptsize $L_1,\;   t \pd_t + \frac{x}{2} \pd_x $   } \\

{\scriptsize  15.}& {\scriptsize   $u_x^2$ } & {\scriptsize $e^u$ } & {\scriptsize $L_0,\;   2t \pd_t + x \pd_x -4 \pd_u$     } \\

{\scriptsize 16.}& {\scriptsize  $u_x^3 + \rho u_x^{-2}$  } & {\scriptsize $\delta$ } & {\scriptsize $L_1,\; \frac{12 t}{11\delta} \pd_t - \frac{6x}{11 \delta  }\pd_x + \left(t^2+  \frac{2 u}{11 \delta  }\right) \pd_u$  }\\

{\scriptsize 17.}& {\scriptsize  $u_x^3 + \rho u_x^{-2}$  } & {\scriptsize $0$ } & {\scriptsize $L_1,\; 2t \pd_t +x \pd_x +  \frac{1}{3} u  \pd_u$   }\\\hline \\[1mm]

\end{tabular}

\end{minipage}
\parbox[t]{0.9\textwidth}{  \protect \footnotesize In the expression of the symmetry generator for Class 8, $k_1, k_2,$ and $k_3$ are arbitrary constants while the function $F$ is a solution of the linear equation $A u_{xxxx} + 2 \beta u_{xx} - u_{tt} + \delta_0 u - \delta k_3=0.$ In the table, the constant parameters $\beta$ and  $\delta$ are arbitrary, $\rho$ and $\theta$ are nonzero, $\alpha \neq 0,1$ and $\beta_0, \delta_0 \in \set{0,1},$ with $\beta \, \beta_0=0.$ }
\end{table}

It should be noted that any value of $S= S(u_x)$ extends the principal algebra by two dimensions when $Z\sim \delta$ or $Z\sim u,$ where $\delta$ is a scalar. In other words \eqref{e:main} contains two non equivalent families of equations indexed by the arbitrary function $S,$ each of which has a symmetry algebra extending the principal algebra by two dimensions. \par

\section{Exact solutions}
\label{s:soln}
In this section we derive exact solutions of a model equation from the family \eqref{e:main} of equations. Namely, we derive these solutions for the fifth equation  from the classification results in  Table \ref{tb:1}, which is given by
\begin{equation}\label{e:5th}
u_{tt} -(2 \beta + 3 u_x) u_{xx} - A u_{xxxx}=0,
\end{equation}
where $\beta$ and $A$ are nonzero real constants. We note that the scaling transformation $t= \sqrt{A} y, x= \sqrt{A} z, u = \sqrt{A} w$ maps \eqref{e:5th}  into exactly the same equation  in which the parameter $A$ is now however reduced to $1.$ In other words, one may assume without loss of generality as we shall   do that $A=1$ in \eqref{e:5th}.

\subsection{Travelling waves through order reductions}
\label{ss:solitary}

Given that the principal algebra of the main equation \eqref{e:main} has generators $\pd_t$ and $\pd_x,$  the family of equations \eqref{e:main} is invariant under both space and time translations. Consequently each equation from this family may in principle have travelling waves solutions, an this holds in particular for \eqref{e:5th}.\par

A remarkable fact about Equation \eqref{e:main} is that it is a Lagrange equation, that is the Euler-Lagrange equation $E_u(\mathscr{L})=0$ of some variational problem with Lagrangian $\mathscr{L},$ where $E_u$ is the Euler operator with respect to $u.$ This is in particular the case for all classes of equations listed  in Table \ref{tb:1}. This property, combined with the fact that the $(1+1)$-dimensional fourth order Equation \eqref{e:main} possesses travelling wave solutions imply that under certain conditions certain solutions  of  \eqref{e:1ode5} may be found by quadrature from those of a reduced second order ordinary differential equation ({\sc ode}).\par

We investigate this property of double order reduction for Lagrange equations in the case of the stated equation \eqref{e:5th}.
 A Lagrangian of \eqref{e:5th} of the same order as the equation is readily found by the standard methods \cite{olv93} to be given by
\begin{equation}\label{e:L0}
 \mathscr{L}_0=\frac{1}{2} \left[  u_{tt} - 2 (\beta + u_x) u_{xx} -  u_{xxxx}  \right].
\end{equation}

Discarding null Lagrangians of order higher than the second from the above expression \eqref{e:L0} gives the second order Lagrangian
\begin{equation}\label{e:lag5th}
\mathscr{L}_1  =  \frac{1}{2} \left[-  u_{xx}^2  +  u_{tt} - 2 (\beta + u_x) u_{xx}  \right]
\end{equation}
for \eqref{e:5th}.  Let $\sigma > 0.$ Then the invariance of  \eqref{e:5th} under the variational symmetry operator
\begin{equation}\label{e:red1}
\bv= \pd_t + \sigma \pd_x \quad \text{ yields } \quad z= x- \sigma t \text{ and } S= u
\end{equation}
as  similarity variables, in terms of which   \eqref{e:5th} reduces to the fourth order ordinary differential equation

\begin{equation}\label{e:1ode5}
(-2 \beta + \sigma^2- 3 S_z) S_{zz}-  S_{zzzz} = 0.
\end{equation}
The latter equation turns out to also be Lagrangian,  and has a second order Lagrangian given by
\begin{equation}\label{e:lag1ode5}
  \mathscr{L}_2 = -S (\beta - \frac{\sigma^2}{2} + S_z)S_{zz} -\frac{1}{2}  S_{zz}^2.
\end{equation}
The hodograph transformation $z= T(\theta),\; S= \theta $ followed by the change of depend variable $T(\theta) = \int V(\theta) d \theta$ reduce $\mathscr{L}_2$ to the first order Lagrangian
\begin{equation} \label{e:lag2ode5}
  \mathscr{L}_3 = \frac{\theta\, V'}{ 8 V^4} +  \frac{\theta (2 \beta -\sigma^2) V'}{ 8 V^3} -  \frac{   V'^2}{ 32\, V^6}.
\end{equation}
By a well-known result of \cite{olv93} (see also  \cite{gungor2,gungor3,gungor4,jcn_Nody1}), every solution of \eqref{e:1ode5} can be found by quadrature from the solutions of the second order {\sc ode}
\begin{equation}\label{e:2ode5}
  E_V(\mathscr{L}_3  - \lambda\, V) \equiv -2 V^3 - 16\, \lambda\, V^7 + V^4 (-4 \beta + 2 \sigma^2)- 3  V'^2 +  V V'' =0,
\end{equation}
for some constant parameter $\lambda,$  where $E_V$ is the Euler operator with respect to the variable $V= V(\theta).$  However, solving  \eqref{e:2ode5} directly is a tedious task. Therefore, we first make the simplifying assumption $\beta = \sigma^2 /2$ in  \eqref{e:2ode5}.   Next, setting $V= 1/\sqrt{w}$ reduces \eqref{e:2ode5} to the slightly simpler equation
\begin{equation}\label{e:3ode5}
  -32 \lambda - 4 w^2 - 3 w^{3/2} w''=0.
\end{equation}
The later equation can be reduced to the integral form
\begin{equation}\label{e:4ode5}
\left(\int_1^{w(x)} \frac{1}{\sqrt{-8 \left(\frac{2 v^{3/2}}{3}-\frac{16 \lambda }{\sqrt{v}}\right)+a_1}} \, dv\right)^2=\left(x+a_2\right)^2,
\end{equation}
for some constants of integration $a_1$ and $a_2.$ Hopefully, specific  solutions $w$ of \eqref{e:4ode5} depending explicitly on $\lambda$ can be found. For each such solution $w= w(\theta),$ one has
\begin{equation}\label{e:1sol5}
 z= T(\theta)= \int \frac{1}{\sqrt{w(\theta)}} d \theta,\quad \text{ and } S= \theta= T^{-1} (z).
\end{equation}
The solution $u$ to \eqref{e:5th} is then obtained by replacing in \eqref{e:1sol5} $S$ by $u$ and $z$ by $x - \sigma t.$ That is
\begin{equation}\label{e:2sol5}
u= T^{-1} (x- \sigma t)
\end{equation}
is the travelling wave solution sought for  \eqref{e:5th}.\par

The expression for $u$ in \eqref{e:2sol5} will yield a  class of solutions  of \eqref{e:main} provided that \eqref{e:4ode5} can be solved explicitly.\par


 Similarly the invariance of  \eqref{e:5th} under the variational symmetry operator
\begin{equation}\label{e:red2}
\bv= \pd_t  \quad \text{ yields } \quad z= x \text{ and } S= u
\end{equation}
as  similarity variables.  In terms of these new variables  \eqref{e:5th} reduces to the new fourth order ordinary differential equation

\begin{equation}\label{e:Xode5}
(-2 \beta - 3 S_z) S_{zz}-  S_{zzzz} = 0
\end{equation}
which is also a Lagrange equation, and has Lagrangian
\begin{equation}\label{e:lagXode5}
  \mathscr{L}_4 = -S S_{zz} (\beta  + S_z) -\frac{1}{2}  S_{zz}^2.
\end{equation}
Using exactly the same double order reduction procedure and same change of variables through which \eqref{e:1ode5} was reduced to the second order counterpart \eqref{e:2ode5}, the corresponding second order reduced equation for \eqref{e:Xode5} is found to be
\begin{equation}\label{e:X2ode5}
  2 2 V^3 - 2 \beta V^4 - 2 \lambda V^6 - 5 V'^2 + 2 V V''=0.
\end{equation}
This time we can more easily find a particular solution $V$ of \eqref{e:X2ode5} for $\lambda=0$, given by
\[
V= - \frac{4}{ a_1^2 + 2 a_1 y + y^2 + 8 \beta}
\]
where  $a_1$ is a constant of integration. The transformation \eqref{e:1sol5} then gives the corresponding solutions $S$ for \eqref{e:Xode5} and $u$ of \eqref{e:5th}  as
\begin{equation}\label{e:solnRedX}
  S=  - 2 \sqrt{2} \sqrt{\beta} \tan\left( \frac{z \sqrt{\beta}+ a_1}{\sqrt{2}}  \right)\quad  \text{ and } u= a_1 - 2 \sqrt{2} \sqrt{\beta} \tan\left( \frac{x \sqrt{\beta}}{\sqrt{2}}  \right),
\end{equation}
and we shall assume that $a_1=0$ in \eqref{e:solnRedX}.  Although the latter solution $u$ expressed in terms of the sole variable $x$ is only a degenerated travelling wave solution, it shows more clearly the validity of the double reduction procedure outlined above.\par
In the same way, the reduction of \eqref{e:5th} under the variational symmetry operator $\bv= \pd_x$ yields the degenerated second order Lagrange equation in $S= S(t)$ given by
\begin{equation}\label{e:Tode5}
  S_{tt}= 0
\end{equation}
with trivial solutions.\par
It is worthwhile  to mention that a  remarkable fact which has appeared up to this point is that the similarity reduction of the Lagrange equation \eqref{e:5th} by the variational symmetry operators $\pd_t, \pd_x$ and $\pd_t + \sigma \pd_x$ has always given rise to  a reduced equation which is also of Lagrange type. This turns out to be a general fact in the actual case of \eqref{e:5th} as each generator of its variational symmetry group actually reduces it to another Lagrange equation. The only verification now needed to be done to confirm this fact is to consider a variational symmetry operator of the form $\bv= \tau \pd_t + \sigma \pd_x$ with $\tau\, \sigma \neq 0.$ It is indeed found that in terms of the resulting similarity variables $z= \tau x  - \sigma t $ and $S=u,$ \eqref{e:5th} reduces to
\begin{equation}\label{e:XTode5}
  (\sigma^2 - 2 \beta\, \tau^2 - 3 \tau^3\, S_z) S_{zz} - \tau^4\, S_{zzzz}=0,
\end{equation}
which is also a Lagrange equation as its Frechet derivative is self-adjoint.   By virtue of the two vector fields $\bv$ used for \eqref{e:main} to \eqref{e:1ode5} and to reducing \eqref{e:XTode5} respectively, the later equation reduces indeed to \eqref{e:1ode5} by the substitution $\tau=1.$  \par

\subsection{Solutions by the generalized tanh method}
\label{ss:solitons}
One can however find more directly solutions of a pre-specified form of  \eqref{e:1ode5}     and hence of \eqref{e:main}   by means of  test based methods often referred to as the tanh methods. In these methods \cite{kudrya,tala}, the solution $S$ of \eqref{e:1ode5} is sought in the form
\begin{equation}\label{e:Tanh1}
  S= \sum_{i=0}^m a_i \left[ g(z; b) \right]^i,
\end{equation}
where $g=g(z; b)$ is a  given known function,  $b=\set{ b_0, \dots, b_p}$ and the parameters $a_i$ and $b_j$ for $i=0, \dots, m$ and $j=0,\dots, p$ are scalars to be found. The integer $m$ is found by inserting the expression of $S$ given by \eqref{e:Tanh1} into the equation and then equating all maximal powers of $g.$ Once $m$ has been found the parameters $a_i$ and $b_j$ are found by solving a system of algebraic equations.  Moreover, given that the solution set of \eqref{e:1ode5}  is invariant under constant shifts, in the search for $S$ in \eqref{e:Tanh1} we may assume that $a_0=0.$\par

      When the function $g=g(z;b)$ in \eqref{e:Tanh1} is given either by $g= a_1 \tanh(b_1 z)$ or $g= a_1/(1+ e^{b_1 z}),$ it is found that $m=1.$
      In the case of the latter function $g$ expressed in terms of the exponential function, the corresponding function $S$ in \eqref{e:Tanh1} is referred to as the exponential rational function \cite{tala}. For these two types of functions $g,$ the corresponding solutions $u= u(t,x)$ of \eqref{e:main} are given  by

\begin{subequations} \label{e:soltanhEq1}
\begin{align}
u &= a_0 \pm 2 \sqrt{- 2 \beta + \sigma^2 }   \tanh \left( \pm \frac{\sqrt{- 2 \beta + \sigma^2 }}{ 2 } (x- \sigma t)  \right)  \label{e:solitV1} \\
\intertext{ and }
u &= a_0  \pm  \frac{ 4 b_0 \sqrt{\sigma^2 - 2 \beta }}{  \exp  \left(\mp   \sqrt{\sigma^2 - 2 \beta } (x- \sigma t)  \right)    + b_0 }  \label{e:solitV2}
\end{align}
\end{subequations}
respectively, where here and in the sequel $a_0$ and $b_0$ are arbitrary constants. The  solutions in \eqref{e:solitV1} are in fact  $1$-solitons \cite{jawad}. Other explicit solutions of \eqref{e:5th} can be found by applying these  test methods.  For example, letting the function $g$ be of the same types  as those used in the preceding case above for the search of solutions of \eqref{e:1ode5} yields solutions of \eqref{e:XTode5} given by
\begin{subequations}\label{e:soltanhEq2}
\begin{align}
S_1(z) &= a_0 \pm \frac{2 \sqrt{\sigma^2 - 2 \beta \tau^2}}{\tau}\tanh \left[ \pm    \frac{ \sqrt{\sigma^2 - 2 \beta \tau^2}}{ 2 \tau^2} z    \right]  \\
\intertext{and}
S_2(z) &= a_0 + \frac{\pm  4 b_0 \sqrt{\sigma^2 - 2 \beta \tau^2}}{ \left[ \exp  \left( \frac{\mp   \sqrt{\sigma^2 - 2 \beta \tau^2} }{\tau^2} \;  z  \right)    + b_0 \right] \tau}.
\end{align}
\end{subequations}
The corresponding solutions of \eqref{e:5th} are then obtained  as $u=S_1(z)$ and $u=S_2(z)$ with  the substitution $z= \tau x- \sigma t.$

\subsection{Multi-soliton solutions by the Hirota method}
We shall now derive families of one- and two-soliton solutions of \eqref{e:5th} using Hirota method \cite{hirot1,hirot2}. One of the crucial steps for obtaining such solutions is to find the appropriate transformation of the dependent variable $u$ in \eqref{e:5th} in terms of one or two auxiliary variables which we may denote here by $\zeta= \zeta(t,x)$ and $\psi= \psi(t,x).$  Indeed, the transformed equation must be expressible as a system of polynomial equations in the so-called Hirota bilinear operators for the unknown functions $\zeta$ and $\psi.$ Ideally, one attempts to express a scalar equation as a quadratic form in terms of a polynomial in the Hirota operators. These operators are pseudo-differential operators $D$ given for a given independent variable $x$ and a pair of functions $\zeta= \zeta(x)$ and $\psi=\psi(x)$ by
\begin{equation}\label{e:H1}
D_x(\zeta, \psi) \equiv D_x\,  \zeta\cdot \psi = \left( \pd_x -\pd_{x'}  \right) \zeta(x) \psi(x') \bigg\vert_{x'=x},
\end{equation}
where $\pd_a$ stands for the partial differential operator $\pd / \pd a$ for any variable $a.$   The expression in \eqref{e:H1} remains clearly unchanged if the functions $\zeta$ and $\psi$ also have independent variables other than $x.$ Higher-order versions of the $D$ operator and its extension to functions of several variables can similarly be defined. In particular for functions $\zeta$ and $\psi$ depending on $x$ and $t$ we have for all nonnegative integers $r$ and $s$

\begin{equation}\label{e:H2}
D_t^r D_x^s \, \zeta \cdot\psi = \left( \pd_x -\pd_{x'}  \right)^s  \left( \pd_t -\pd_{t'}  \right)^r  \zeta(t,x)\psi(t', x')  \bigg\vert_{x'=x, t'=t}.
\end{equation}

The fundamental difference between the Hirota operator $D_x$ and the usual partial differential operator $\pd_x$ is the minus sign that occurs in \eqref{e:H1} in place of the plus sign if $D_x$ is replaced with $\pd_x$ in that equation. In other words, the Hirota $D$ operator does not satisfy the Leibniz rule on the product of two functions as the operator $\pd_x$ does. In particular, although it is a bilinear operator it is not symmetric, but it naturally also acts linearly on the space of functions. Several important properties of the Hirota $D$ operator can be found in \cite{hirot1,hirot2}. See also \cite{hirot3,wazwaz1} and the references therein. For brevity, knowledge of these properties will be assumed in the sequel. \par

More succinctly, in the actual case of \eqref{e:5th} we seek a change of variable of the form $u= \vartheta (f)$ that transforms \eqref{e:5th} into a quadratic form $P(D) f\cdot f,$ where $P(D)$ is a polynomial in the D operators with constant coefficients. To this end, we set
\begin{equation}\label{e:transH1}
u= 2 \alpha \pd_x \log(f), \quad \text{ where } f=f(t,x).
\end{equation}
Then, under the latter transformation \eqref{e:transH1} with $\alpha=2$, the transformed version of  \eqref{e:5th} after integration with respect to the variable $x$ can be put into the form
\begin{equation}\label{e:transH2}
(- D_t^2 + 2 \beta D_x^2 + D_x^4 ) f\cdot f=0.
\end{equation}
The LHS in the above equation is indeed of the form $P(D) f \cdot f,$ where $P$ is the polynomial  given by $P(x_1, x_2)= -x_1^2 + 2 \beta x_2^2 + x_2^4.$

It is known \cite{hirot1,hirot2} that one-soliton solutions of \eqref{e:5th} can be found by seeking solutions $f$ of \eqref{e:transH2} of the form

\begin{subequations} \label{e:1solit}
\begin{align}
f &= 1+ \exp(\xi),\quad \text{ where } \xi= \sigma x - \omega t\\
\intertext{and}
\omega^2 &= 2 \beta \sigma^2 + \sigma^4, \qquad \text{that is}\quad  P(\omega, \sigma) =0,
\end{align}
\end{subequations}
for some positive constants $\omega$ and $\sigma.$ The latter equation $P(\omega, \sigma)=0$ in \eqref{e:1solit} is called the dispersion relation for \eqref{e:5th}. On the other hand, $\sigma$ is called the wave-number and $\omega$ is the frequency, so that the velocity of the wave crest is $\omega/ \sigma.$  It turns out that the function $f$ given by \eqref{e:1solit} does satisfy  \eqref{e:transH2}. Consequently, substituting \eqref{e:1solit} into \eqref{e:transH1} with $\alpha=2$ yields
\begin{subequations}\label{e:1solitFnal}
\begin{align}
  u &= \frac{4\sigma \exp(\xi)}{1+\exp(\xi)}  \\
    &= 2 \sigma \exp(\xi/2) \sech(\xi/2),
\end{align}
\end{subequations}
and the expression of $u$ in \eqref{e:1solitFnal} represents therefore a one-parameter family of one-soliton solutions of \eqref{e:5th}.\par

Similarly,  two-soliton solutions of \eqref{e:5th} may be found by seeking the solution $f$ of \eqref{e:transH2} in the form
\begin{subequations}\label{e:2solit}
\begin{align}
f  &= 1+ \exp(\xi_1) + \exp(\xi_2) +  \Phi \exp(\xi_1+ \xi_2), \\
\intertext{ where}
\xi_1 &= - \omega_1 t + \sigma_1 x,\quad \xi_2= - \omega_2 t + \sigma_2 x,\\
\omega_1 &= \sqrt{2 \beta \sigma_1^2 + \sigma_1^4},\quad \omega_2 = \sqrt{2 \beta \sigma_2^2 + \sigma_2^4},\quad \Phi= -  \frac{P(\omega_2 -\omega_1, \sigma_1-\sigma_2)}{P(\omega_1+\omega_2, \sigma_1+\sigma_2)},
\end{align}
\end{subequations}
and for some constant parameters $\omega_i$ and  $\sigma_i,$ for $i=1,2.$ The function $f$ thus defined in \eqref{e:2solit} actually also turns out to be a solution of \eqref{e:transH2}.  The substitution of the expression of $f$ from \eqref{e:2solit} into \eqref{e:transH1} therefore yields a solution of \eqref{e:5th} given by
\begin{equation}\label{e:2solitF0}
  u = \frac{4 \left(e^{x \sigma_1 +t  \omega_2 } \sigma_1 +e^{x \sigma_2 +t  \omega_1 } \sigma_2 \right)+4
e^{x ( \sigma_1 +\sigma_2)} (\sigma_1+\sigma_2) \Phi }{e^{x \sigma_2 +t  \omega_1 }
+e^{t (\omega_1 + \omega_2 )}+e^{x \sigma_1 +t  \omega_2 }+e^{x (\sigma_1+\sigma_2)} \Phi }. \\
\end{equation}

The latter solution $u$ therefore represents a   two-soliton solution of \eqref{e:5th}. Moreover, substituting $\omega_1, \omega_2,$ and $\Phi$ in \eqref{e:2solitF0} in terms of their expressions in \eqref{e:2solit} gives the more explicit solution of \eqref{e:5th} of the form

{\small
\begin{align}\label{e:2solitFn}
 &u =  \frac{4\left( e^{x \sigma_1+ t \sqrt{2 \beta \sigma_2^2 + \sigma_2^4}} \  R_1 \sigma_1 + e^{x \sigma_2+ t \sqrt{2 \beta \sigma_1^2 + \sigma_1^4}} \  R_1 \sigma_2  -   e^{x (\sigma_1+\sigma_2)}R_2 (\sigma_1+\sigma_2) \right)     } {\left( \left(e^{x \sigma_2+ t \sqrt{2 \beta \sigma_1^2 + \sigma_1^4}} +    e^{ t (\sqrt{2 \beta \sigma_1^2 + \sigma_1^4} + \sqrt{2 \beta \sigma_2^2 + \sigma_2^4})}  +
  e^{x \sigma_1+ t \sqrt{2 \beta \sigma_2^2 + \sigma_2^4}}\right)  R_1 - e^{x (\sigma_1+ \sigma_2)} R_2 \right)                                       }   \\
  \intertext{where}
 &R_1= \sigma_1 \sigma_2 \left(2(\beta+ \sigma_1^2) + 3 \sigma_1 \sigma_2+ 2 \sigma_2^2\right) - \sqrt{2 \beta \sigma_1^2 + \sigma_1^4 }   \sqrt{2 \beta \sigma_2^2 + \sigma_2^4 } \\
  &R_2 = - \sigma_1 \sigma_2 \left(2(\beta+ \sigma_1^2) - 3 \sigma_1 \sigma_2+ 2 \sigma_2^2\right) + \sqrt{2 \beta \sigma_1^2 + \sigma_1^4 }   \sqrt{2 \beta \sigma_2^2 + \sigma_2^4 }
\end{align}
}

Letting $\beta = 4 \sigma_1^2$ in \eqref{e:2solitFn} and setting  $\sigma_2= 8 \sigma_1$ yields the particular solution $u_p$ of \eqref{e:5th}, whose expression takes the form
  \begin{align}\label{e:solnUp}
u_p=& \frac{12 \sigma_1 e^{\sigma_1 x} \left(k e^{48 \sqrt{2} \sigma_1^2 t}+\left(-57+9 \sqrt{2}\right)
e^{8 \sigma_1 x}+8 k e^{3 \sigma_1^2 t+7 \sigma_1 x}\right)}{3 k e^{\left(3+48 \sqrt{2}\right) \sigma_1^2 t}+\left(-19+3 \sqrt{2}\right)
e^{9 \sigma_1 x}+3 k e^{48 \sqrt{2} \sigma_1^2 t+\sigma_1 x}+3 k e^{3 \sigma_1^2 t+8 \sigma_1 x}}
\end{align}
where $k= \left(-9+\sqrt{2}\right).$  Other particular solutions of \eqref{e:5th} can be obtained from \eqref{e:2solitFn} in a similar manner.\par

\begin{figure}[ht]
\begin{center}
\caption{\label{fg:profileUpV1} \protect \footnotesize
Profiles of the solution $u_p$ in \eqref{e:solnUp} with $\sigma_1=6.3.$  These profiles correspond to  times $t=-5000,\ t=-500,\ t=-50$ and $t=-5.$ The profiles show that the wave keeps a constant shape with a wave amplitude oscillating around $226.8.$
\vspace{2mm}}
\includegraphics[width=\textwidth]{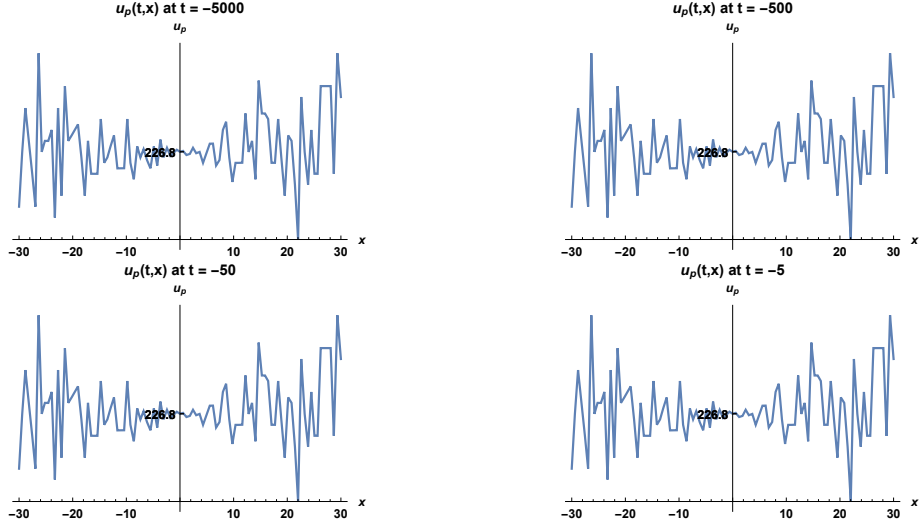}
\end{center}
\end{figure}

\begin{figure}[ht]
\begin{center}
\caption{\label{fg:profileUpV2} \protect \footnotesize
Profiles of the solution $u_p$ in \eqref{e:solnUp} with $\sigma_1=6.3.$ The times corresponding to these profiles are $t=1,\ t=10,$ and $t=5000.$ The profiles show that the wave amplitude eventually vanishes from around time $t=1$ onward. This vanishing is achieved by time $t=10$ and kept up to time $t=5000$ and beyond.
\vspace{2mm}}
\includegraphics[width=\textwidth]{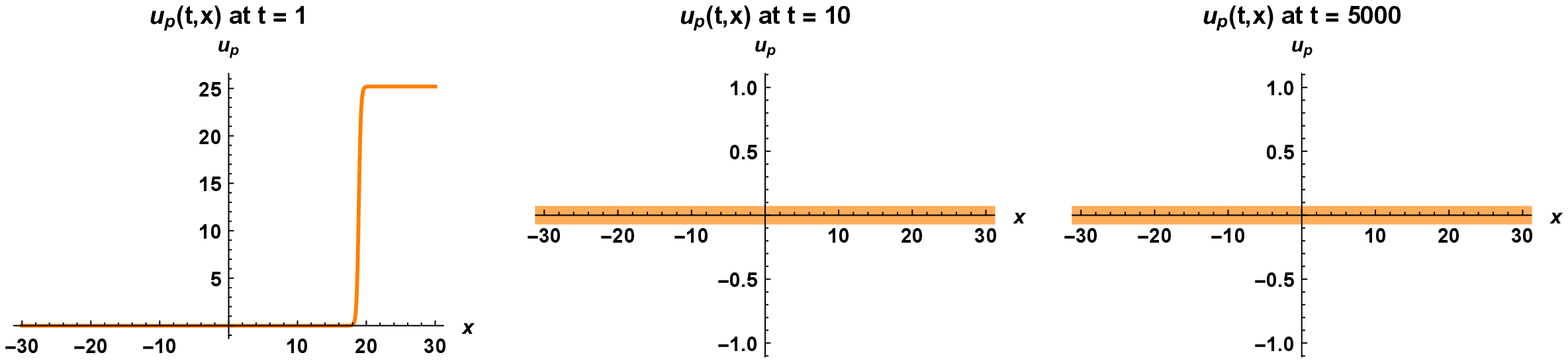}
\end{center}
\end{figure}

Figure \ref{fg:profileUpV1} shows that between time $t\leq -5000$ up to time $t=-5$ the function $u_p$ preserves its shape and this is in fact the case up to around time $t=-2,$ with an amplitude oscillating closely around $226.8$  On the other hand, Figure \ref{fg:profileUpV2} shows that the wave amplitude starts decreasing to zero from around time $t> -2$ and by time $t=10$ the amplitude has vanished, and this zero-amplitude profile is verified in the same figure to remain valid up to time $t=5000.$ The value of $\sigma_1=6.3$ is used for profiling the shape of $u_p$ in Figure \ref{fg:profileUpV1} and Figure \ref{fg:profileUpV2}, although  assigning any value to $\sigma_1$ is only needed for the numerical computation of profiles.

\subsection{Symmetry generated solutions}
\label{ss:symsol}
The very fact that every symmetry transformation of \eqref{e:5th} leaves it invariant implies in particular that every such transformation maps a solution of the equation to another solution. Thus for any given solution, the whole symmetry group generates an $r$-parameter family of new solutions, where $r$ is the dimension of the symmetry group,  which equals five in the case of \eqref{e:5th}. The symmetry transformations themselves are given for each symmetry generator $\bv_j$ by its flow $\Psi_j(\eps, P)= F_j,$ where $P=(t,x,u)$ and $F_j= (y,z, w)$ is the transformation of $P$ under the flow $\Psi_j$ of $\bv_j,$ and $\eps$ is the group parameter.\par

Let us denote the symmetry generators of \eqref{e:5th} by
\begin{align}\label{e:gen5}
  \bv_1 &= \pd_t,\qquad \bv_2= \pd_x, \qquad \bv_3= \pd_u,\qquad \bv_4= t  \pd_u  \\
  \bv_5 &= -\frac{3t}{2 \beta}  \pd_t -    \frac{3x}{4 \beta} \pd_x + \left(x +  \frac{3u}{4 \beta} \right) \pd_u.
\end{align}
The corresponding flows are then given by
\begin{align}\label{e:flow}
F_1 & = (t+ \eps, x, u),\qquad F_2 = (t, x+ \eps, u)\\
F_3& = (t, x, u+ \eps),\qquad  F_4 = (t, x, t \eps + u)\qquad \\
F_5&= \left[ t e^\frac{-3\eps} {2 \beta}, e^\frac{-3\eps} {4 \beta}  x,    -\frac{2}{3} x \beta e^\frac{-3\eps} {4 \beta} +
e^\frac{3\eps} {4 \beta} \left( u+ \frac{2 x \beta}{3} \right) \right].
\end{align}

Each of these group actions transforms indeed any given  solution of \eqref{e:5th} into a solution of the same equation. For instance if $u= H(t,x)$ is a solution of
\eqref{e:5th}, then $F_1$ transforms it into another solution $u=H(t-\eps,x),$ reflecting the invariance of \eqref{e:5th} under time translations. More generally if $k_1, \dots, k_5$ are sufficiently small arbitrary constants, then any  element of the whole symmetry group of \eqref{e:5th} is determined by such constants and transforms a solution $u=H(t,x),$ as composite of the actions of each flow   $F_j = \Psi_j (k_j, P)$ on $u= H(t, x)$ to the new solution
\begin{equation}\label{e:grpSoln}
  u= k_3 + k_4 t + \frac{2}{3}\left( -1 + e^{\frac{3 k5}{2 \beta}} \right) x \beta + e^{\frac{3 k5}{4 \beta}}H\left( t e^{\frac{3 k5}{2 \beta}} -k_1, x e^{\frac{3 k5}{4 \beta}} - k_2  \right).
\end{equation}
This shows amongst others that adding a linear function of $x$ and $t$ to any given solution yields another solution of the equation. For instance, the degenerated solution $u$ in \eqref{e:solnRedX} depending only of the space variable $x$ is transformed under the combined action of $\bv_4$ and $\bv_5$ into a new solution depending on all independent variables in a nontrivial manner, of the form
\begin{equation}\label{e:NewTan}
u= \delta  t + \frac{2}{3} x \beta (\rho^2 -1) - \rho \left(k_2 + 2 \sqrt{2} \sqrt{\beta} \tan\left( \frac{\rho x \sqrt{\beta}}{\sqrt{2}}\right)\right),
\end{equation}
where $\delta$ and $\rho$ are arbitrary parameters with $\rho>0.$ This outlines the fact that  starting even with a trivial solution of the equation, the symmetry transformations can usually generate the needed or sufficiently insightful solution to the problem at hand.\par

Moreover \eqref{e:grpSoln} gives a clear answer to the crucial problem of generation of new soliton solutions. Indeed, if we let $k_4= k_5=0$  in \eqref{e:grpSoln}, then the resulting transformation

\begin{equation}\label{e:grpSoliton}
  u= k_3 +  H\left( t  -k_1, x - k_2  \right)
\end{equation}
preserves soliton solutions, outlining the fact that such solutions are preserved under translations of dependent and independent variables.


\section{Concluding Remarks}
\label{s:conclusion}
In this paper, we considered for group classification a particular family \eqref{e:main} of the general family
\begin{equation}\label{e:gnl1p1Wave}
  u_{tt} = F(u, u_x, u_{xx}, u_{xxx}, u_{xxxx}, \dots)
\end{equation}
of nonlinear wave equations in $(1+1)$ dimensions \cite{anco-kara}, where $F$ is a given function of its arguments.  This study  has revealed amongst others the effectiveness of our method of indeterminates for group classification in the case of equations depending on two arbitrary labeling functions with each depending on its own argument. We've also given here a more formal and clearer description of the method.  Needless to say that the method should work for equations involving an arbitrary number of such labelling functions. Although it is clear that the method can only apply to cases where the determining equations can be expressed in terms of some  indeterminates, the success of the method to a variety of equations including linear ones now leads us to explicitly look at amongst others the main types of equations for which the method could not be applied.\par

By exploiting the Lagrange nature of the subfamily \eqref{e:5th} of the fourth order equation \eqref{e:main}, we have been able to reduce  it to a second order {\sc ode}, first through a similarity reduction to a fourth order {\sc ode} via variational symmetries,  followed by  a double order reduction.  Solitons and other   solutions of the model equation \eqref{e:5th} have also been found by four different methods including Lie group and Hirota methods, and the most general action of the symmetry group on a given solution has been explicitly determined.\par

Our profiling of the particular solution $u_p$ of \eqref{e:5th} given by \eqref{e:solnUp} and depicted in Figures \ref{fg:profileUpV1} and \ref{fg:profileUpV2} shows clearly as expected for solitons that $u_p$ preserves its shape throughout time. In \cite{kumar} however, and quite often in the literature (see the references in \cite{kumar}), the plot of solutions duly referred to as solitons displays different profiles for different fixed times. In reality, very few methods of solutions are expected to produce solitons. Although the Hirota method is one of them, solutions produced by the popular the Lie symmetry methods have \emph{a priori} no connections with solitons. On the other hand, although the solution $u_p$ depicted in Figure  \ref{fg:profileUpV2} eventually vanishes, it is not a breather solution as studied in \cite{scheider} and in some of the references therein on breather solutions. Indeed, the vanishing of $u_p$ in \eqref{e:solnUp} occurs in a fixed Cartesian frame and not in one that moves with the velocity of the wave function. \par

On the other hand, this study has also given rise in particular to an explicit formula for genuinely generating new soliton solutions from any given one, a result of crucial importance seldom found in the current literature. \par

Finally, it is well-known that a variational problem is invariant under its variational symmetry group, and that the latter transforms the Euler-Lagrange equation into an equivalent Euler-Lagrange equation  having as Lagrangian the a multiple of the original one. No similar results seem to be known with regard to the similarity reduction of Euler-Lagrange equations. It has however emerged from this study that the similarity reduction of the model Lagrange equation \eqref{e:5th} by any generator of its variational symmetry group also  yields a Lagrange equation. This is a result whose generalization we wish to investigate in our future research.


\section*{Declarations of interest:} None.
\label{s:declare}


\bibliographystyle{model1-num-names}


\end{document}